\documentclass{amsart}

\usepackage{amsfonts,amssymb,verbatim,amsmath,amsthm,latexsym,textcomp,amscd}
\usepackage{latexsym,amsfonts,amssymb,epsfig,verbatim}
\usepackage{amsmath,amsthm,amssymb,latexsym,graphics,textcomp}
\usepackage{graphicx}
\usepackage{color}
\usepackage{url}

\input xy
\xyoption{all}

\begin{document}

\newtheorem{theorem}{Theorem}[section]
\newtheorem{prop}[theorem]{Proposition}
\newtheorem{lemma}[theorem]{Lemma}
\newtheorem{cor}[theorem]{Corollary}
\newtheorem{definition}[theorem]{Definition}
\newtheorem{conj}[theorem]{Conjecture}
\newtheorem{rmk}[theorem]{Remark}
\newtheorem{claim}[theorem]{Claim}
\newtheorem{defth}[theorem]{Definition-Theorem}

\newcommand{\boundary}{\partial}
\newcommand{\C}{{\mathbb C}}
\newcommand{\integers}{{\mathbb Z}}
\newcommand{\natls}{{\mathbb N}}
\newcommand{\ratls}{{\mathbb Q}}
\newcommand{\bbR}{{\mathbb R}}
\newcommand{\proj}{{\mathbb P}}
\newcommand{\lhp}{{\mathbb L}}
\newcommand{\tube}{{\mathbb T}}
\newcommand{\cusp}{{\mathbb P}}
\newcommand\AAA{{\mathcal A}}
\newcommand\BB{{\mathcal B}}
\newcommand\CC{{\mathcal C}}
\newcommand\DD{{\mathcal D}}
\newcommand\EE{{\mathcal E}}
\newcommand\FF{{\mathcal F}}
\newcommand\GG{{\mathcal G}}
\newcommand\HH{{\mathcal H}}
\newcommand\II{{\mathcal I}}
\newcommand\JJ{{\mathcal J}}
\newcommand\KK{{\mathcal K}}
\newcommand\LL{{\mathcal L}}
\newcommand\MM{{\mathcal M}}
\newcommand\NN{{\mathcal N}}
\newcommand\OO{{\mathcal O}}
\newcommand\PP{{\mathcal P}}
\newcommand\QQ{{\mathcal Q}}
\newcommand\RR{{\mathcal R}}
\newcommand\SSS{{\mathcal S}}
\newcommand\TT{{\mathcal T}}
\newcommand\UU{{\mathcal U}}
\newcommand\VV{{\mathcal V}}
\newcommand\WW{{\mathcal W}}
\newcommand\XX{{\mathcal X}}
\newcommand\YY{{\mathcal Y}}
\newcommand\ZZ{{\mathcal Z}}
\newcommand\CH{{\CC\HH}}
\newcommand\PEY{{\PP\EE\YY}}
\newcommand\MF{{\MM\FF}}
\newcommand\RCT{{{\mathcal R}_{CT}}}
\newcommand\PMF{{\PP\kern-2pt\MM\FF}}
\newcommand\FL{{\FF\LL}}
\newcommand\PML{{\PP\kern-2pt\MM\LL}}
\newcommand\GL{{\GG\LL}}
\newcommand\Pol{{\mathcal P}}
\newcommand\half{{\textstyle{\frac12}}}
\newcommand\Half{{\frac12}}
\newcommand\Mod{\operatorname{Mod}}
\newcommand\Area{\operatorname{Area}}
\newcommand\ep{\epsilon}
\newcommand\hhat{\widehat}
\newcommand\Proj{{\mathbf P}}
\newcommand\U{{\mathbf U}}
 \newcommand\Hyp{{\mathbf H}}
\newcommand\D{{\mathbf D}}
\newcommand\Z{{\mathbb Z}}
\newcommand\R{{\mathbb R}}
\newcommand\Q{{\mathbb Q}}
\newcommand\E{{\mathbb E}}
\newcommand\til{\widetilde}
\newcommand\length{\operatorname{length}}
\newcommand\tr{\operatorname{tr}}
\newcommand\gesim{\succ}
\newcommand\lesim{\prec}
\newcommand\simle{\lesim}
\newcommand\simge{\gesim}
\newcommand{\simmult}{\asymp}
\newcommand{\simadd}{\mathrel{\overset{\text{\tiny $+$}}{\sim}}}
\newcommand{\ssm}{\setminus}
\newcommand{\diam}{\operatorname{diam}}
\newcommand{\pair}[1]{\langle #1\rangle}
\newcommand{\T}{{\mathbf T}}
\newcommand{\inj}{\operatorname{inj}}
\newcommand{\pleat}{\operatorname{\mathbf{pleat}}}
\newcommand{\short}{\operatorname{\mathbf{short}}}
\newcommand{\vertices}{\operatorname{vert}}
\newcommand{\collar}{\operatorname{\mathbf{collar}}}
\newcommand{\bcollar}{\operatorname{\overline{\mathbf{collar}}}}
\newcommand{\I}{{\mathbf I}}
\newcommand{\tprec}{\prec_t}
\newcommand{\fprec}{\prec_f}
\newcommand{\bprec}{\prec_b}
\newcommand{\pprec}{\prec_p}
\newcommand{\ppreceq}{\preceq_p}
\newcommand{\sprec}{\prec_s}
\newcommand{\cpreceq}{\preceq_c}
\newcommand{\cprec}{\prec_c}
\newcommand{\topprec}{\prec_{\rm top}}
\newcommand{\Topprec}{\prec_{\rm TOP}}
\newcommand{\fsub}{\mathrel{\scriptstyle\searrow}}
\newcommand{\bsub}{\mathrel{\scriptstyle\swarrow}}
\newcommand{\fsubd}{\mathrel{{\scriptstyle\searrow}\kern-1ex^d\kern0.5ex}}
\newcommand{\bsubd}{\mathrel{{\scriptstyle\swarrow}\kern-1.6ex^d\kern0.8ex}}
\newcommand{\fsubeq}{\mathrel{\raise-.7ex\hbox{$\overset{\searrow}{=}$}}}
\newcommand{\bsubeq}{\mathrel{\raise-.7ex\hbox{$\overset{\swarrow}{=}$}}}
\newcommand{\tw}{\operatorname{tw}}
\newcommand{\base}{\operatorname{base}}
\newcommand{\trans}{\operatorname{trans}}
\newcommand{\rest}{|_}
\newcommand{\bbar}{\overline}
\newcommand{\UML}{\operatorname{\UU\MM\LL}}
\newcommand{\EL}{\mathcal{EL}}
\newcommand{\tsum}{\sideset{}{'}\sum}
\newcommand{\tsh}[1]{\left\{\kern-.9ex\left\{#1\right\}\kern-.9ex\right\}}
\newcommand{\Tsh}[2]{\tsh{#2}_{#1}}
\newcommand{\qeq}{\mathrel{\approx}}
\newcommand{\Qeq}[1]{\mathrel{\approx_{#1}}}
\newcommand{\qle}{\lesssim}
\newcommand{\Qle}[1]{\mathrel{\lesssim_{#1}}}
\newcommand{\simp}{\operatorname{simp}}
\newcommand{\vsucc}{\operatorname{succ}}
\newcommand{\vpred}{\operatorname{pred}}
\newcommand\fhalf[1]{\overrightarrow {#1}}
\newcommand\bhalf[1]{\overleftarrow {#1}}
\newcommand\sleft{_{\text{left}}}
\newcommand\sright{_{\text{right}}}
\newcommand\sbtop{_{\text{top}}}
\newcommand\sbot{_{\text{bot}}}
\newcommand\sll{_{\mathbf l}}
\newcommand\srr{_{\mathbf r}}
\newcommand\geod{\operatorname{\mathbf g}}
\newcommand\mtorus[1]{\boundary U(#1)}
\newcommand\A{\mathbf A}
\newcommand\Aleft[1]{\A\sleft(#1)}
\newcommand\Aright[1]{\A\sright(#1)}
\newcommand\Atop[1]{\A\sbtop(#1)}
\newcommand\Abot[1]{\A\sbot(#1)}
\newcommand\boundvert{{\boundary_{||}}}
\newcommand\storus[1]{U(#1)}
\newcommand\Momega{\omega_M}
\newcommand\nomega{\omega_\nu}
\newcommand\twist{\operatorname{tw}}
\newcommand\modl{M_\nu}
\newcommand\MT{{\mathbb T}}
\newcommand\Teich{{\mathcal T}}
\renewcommand{\Re}{\operatorname{Re}}
\renewcommand{\Im}{\operatorname{Im}}

\title{Dynamics of $L^p$ multipliers on harmonic manifolds}

\author{Kingshook Biswas}
\address{Indian Statistical Institute, Kolkata, India. Email: kingshook@isical.ac.in}
\author{Rudra P. Sarkar}
\address{Indian Statistical Institute, Kolkata, India. Email: rudra@isical.ac.in}

\begin{abstract} Let $X$ be a complete, simply connected harmonic manifold with sectional
curvatures $K$ satisfying $K \leq -1$. In \cite{biswas6}, a Fourier transform was defined for
functions on $X$, and a Fourier inversion formula and Plancherel theorem were proved. We use the
Fourier transform to investigate the dynamics on $L^p(X)$ for $p > 2$ of certain bounded linear
operators $T : L^p(X) \to L^p(X)$ which we call "$L^p$-multipliers" in accordance with standard
terminology. These operators are required to preserve the subspace of $L^p$ radial functions.
A notion of convolution with radial functions was defined in \cite{biswas6},
and these operators are also required to be compatible with convolution in the sense that
$$
T\phi * \psi = \phi * T\psi
$$
for all radial $C^{\infty}_c$-functions $\phi, \psi$. They are also required to be
compatible with translation of radial functions.
Examples of $L^p$-multipliers are given by the
operator of convolution with an $L^1$ radial function, or more generally convolution with a finite
radial measure. In particular elements of the heat semigroup $e^{t\Delta}$ act as multipliers.
Given $2 < p < \infty$, we show that for any $L^p$-multiplier $T$ which is not a scalar multiple of the identity, there is an
open set of values of $\nu \in \C$ for which the operator $\frac{1}{\nu} T$ is chaotic on $L^p(X)$ in the sense of Devaney,
i.e. topologically transitive and with periodic points dense. Moreover such operators are
topologically mixing. We also show that there is a constant $c_p > 0$ such that for any $c \in \C$ with $\Re c > c_p$,
the action of the shifted heat semigroup $e^{ct} e^{t\Delta}$ on $L^p(X)$ is chaotic. These results generalize the
corresponding results for rank one symmetric spaces of noncompact type and negatively curved harmonic $NA$ groups
(or Damek-Ricci spaces).
\end{abstract}

\bigskip

\maketitle

\tableofcontents

\section{Introduction}

\medskip

The study of chaos in linear dynamics originated in the work of Godefroy and Shapiro \cite{godefroyshapiro}.
The dynamics of a linear operator $T$ on a Frechet space $X$ is said to be {\it chaotic} (in the sense
of Devaney) if $T$ is {\it hypercyclic} (i.e. has a dense orbit, equivalently is topologically transitive),
and has a dense set of periodic points.
There is now an extensive literature on chaotic and hypercyclic operators, of which a summary may be found
in the books \cite{bayartmatheron}, \cite{erdmannmanguillot}.

\medskip

In a geometric context, linear chaos has been investigated for the heat semigroup $e^{t\Delta}$ acting on
the Lebesgue spaces $L^p(X)$, for certain complete Riemannian manifolds $X$ (where $\Delta = div \ grad$ is the Laplace-Beltrami
operator on $X$). Ji and Weber considered finite volume locally symmetric spaces of rank one in \cite{jiweber1},
where they showed that for $p \in (1,2)$ there is a constant $c_p \in \R$ such that for $c > c_p$ the shifted
semigroup $e^{t(\Delta + c)}$ is {\it subspace chaotic} on $L^p(X)$, i.e. there is a closed, invariant subspace
such that the semigroup restricted to the subspace is chaotic. In \cite{jiweber2}, Ji and Weber investigated the case of
symmetric spaces of noncompact type, and showed that in this setting
for $p \in (2, \infty)$ there is a constant $c_p \in \R$ such that for $c > c_p$ the shifted
semigroup $e^{t(\Delta + c)}$ is subspace chaotic on $L^p(X)$. In \cite{sarkar}, Sarkar improved
on the result of Ji and Weber for rank one symmetric spaces, by showing that for the {\it Damek-Ricci spaces}
(these are certain solvable Lie groups equipped with a left-invariant metric, which include as a particular
case rank one symmetric spaces of noncompact type \cite{damekricci1}),
for $p \in (2, \infty)$ there is a constant $c_p \in \R$ such that for $c > c_p$ the shifted
semigroup $e^{t(\Delta + c)}$ is chaotic on $L^p(X)$, and not just subspace chaotic. Sarkar and
Pramanik later showed that the same result also holds for higher rank symmetric spaces of noncompact type \cite{sarkarpramanik}.
Ji and Weber also extended their results for locally symmetric spaces to the case of higher rank in \cite{jiweber3}.
Finally, in \cite{sarkarray}, Sarkar and Ray generalized the results on chaotic dynamics of the heat semigroup to the
case of more general operators on symmetric spaces of noncompact type known as {\it Fourier multipliers} (these include
as a particular case the operators $e^{t\Delta}$), showing that
for $p \in (2, \infty)$, for any such operator $T$ on $L^p(X)$ which is not a
scalar multiple of the identity, there is a $z \in \C$ such that the operator $zT$ is chaotic.

\medskip

The aim of the present article is to generalize this last result to the case of a class of Riemannian
manifolds known as {\it harmonic manifolds}. These include the rank one symmetric spaces and Damek-Ricci spaces
as particular examples. A Riemannian manifold $X$ is said to be {\it harmonic} if for any $x \in X$, sufficiently
small geodesic spheres centered at $x$ have constant mean curvature depending only on the radius of the sphere.
Harmonic manifolds may be characterized in various equivalent ways, one characterization being that
harmonic functions on the manifold satisfy the mean-value property with respect to geodesic spheres. The Lichnerowicz
conjecture asserts that harmonic manifolds are either flat or locally symmetric of rank one. The conjecture holds
in dimension less than or equal to 5 (\cite{copsonruse}, \cite{walker1}, \cite{nikolayevsky}) and for
compact simply connected harmonic manifolds \cite{szabo}, though it is false in general, with the Damek-Ricci
spaces giving a family of counterexamples \cite{damekricci1}. Heber showed however that the only complete,
simply connected, homogeneous harmonic manifolds are the Euclidean spaces, rank one symmetric spaces, and the Damek-Ricci
spaces \cite{heber}. For a survey of results on general noncompact harmonic manifolds we refer
to \cite{knieperpeyerimhoff1}.

\medskip

In \cite{biswas6}, a study of harmonic analysis on noncompact harmonic manifolds in terms of eigenfunctions
of the Laplace-Beltrami operator $\Delta$ was initiated. We briefly describe the results from \cite{biswas6}
which we will be needing:

\medskip

Let $X$ be a complete, simply connected, harmonic manifold with sectional curvatures $K$ satisfying $K \leq -1$.
Then $X$ is a CAT(-1) space, and can be compactified by adjoining a {\it boundary at infinity} $\partial X$, given by
equivalence classes of geodesic rays $\gamma : [0, \infty) \to X$ in $X$. The {\it Busemann cocycle} $B : X \times X \times \partial X \to \R$
is defined by
$$
B(x, y, \xi) := \lim_{z \in X, z \to \xi} (d(x, z) - d(y, z))
$$
Given $x \in X$ and $\xi \in \partial X$, the {\it Busemann function} at $\xi$ based at $x$ is defined by
$B_{\xi, x}(y) := B(y, x, \xi)$. The Busemann functions $B_{\xi, x}$ are $C^2$ convex functions, and their level
sets are called {\it horospheres} based at $\xi$. As $X$ is harmonic, $X$ is also {\it asymptotically harmonic},
i.e. all horospheres have constant mean curvature, so there is a constant $h \in \R$ such that $\Delta B_{\xi, x} \equiv h$
for all $x \in X, \xi \in \partial X$. Since $X$ is negatively curved, in fact $h > 0$. We let
$$
\rho = \frac{1}{2}h
$$
Then in \cite{biswas6} it is shown that for all $\lambda \in \C$, the function $e^{(i\lambda - \rho)B_{\xi, x}}$ is an
eigenfunction of $\Delta$ with eigenvalue $-(\lambda^2 + \rho^2)$. Given $f \in C^{\infty}_c(X)$ and $x \in X$,
the {\it Fourier transform} of $f$ based at $x$ is the function $\tilde{f}^x$ on $\C \times \partial X$ defined by
$$
\tilde{f}^x(\lambda, \xi) := \int_X f(y) e^{(-i\lambda - \rho)B_{\xi, x}(y)} dvol(y)
$$

\medskip

Given $x \in X$, a function $f$ on $X$ is said to be {\it radial around $x$} if $f$ is constant on
geodesic spheres centered at $x$. In \cite{biswas6} it is shown that for any $\lambda \in \C$, there is a
unique eigenfunction $\phi_{\lambda, x}$ of $\Delta$ for the eigenvalue $-(\lambda^2 + \rho^2)$ which is
radial around $x$ and satisfies $\phi_{\lambda, x}(x) = 1$. For $p > 2$, the functions $\phi_{\lambda, x}$ are in
$L^p(X)$ for $\lambda$ in the strip
$$
S_p := \{ \lambda \in \C | |\Im \lambda| < (1 - 2/p)\rho \}
$$
Let $1 \leq q < 2$ be such that $1/p + 1/q = 1$. The {\it spherical Fourier transform} based at $x$ of a function $f \in L^q(X)$
is the function $\hat{f}^x$ on $\R$ defined by
$$
\hat{f}^x(\lambda) := \int_X f(y) \phi_{\lambda, x}(y) dvol(y)
$$
The spherical Fourier transform $\hat{f}^x$ extends to a holomorphic function on the strip $S_p$.

\medskip

When $X$ is a rank one symmetric space of noncompact type, an {\it $L^p$-multiplier} is a bounded operator $T : L^p(X) \to L^p(X)$
which is translation invariant. Examples of $L^p$-multipliers are given by convolution on the right with radial $L^1$-functions, or more
generally convolution on the right with finite radial measures. For a general harmonic manifold as in our case, a notion of convolution of functions
with radial functions was defined in \cite{biswas6} as follows: given a function $g$ radial around a point $x$ and another point
$y \in X$, the $y$-translate of $g$ is the function $\tau_y g$ radial around $y$ defined by requiring that the value of $\tau_y g$
on a geodesic sphere of radius $r$ around $y$ equals the value of $g$ on the geodesic sphere of radius $r$ around $x$. The convolution
of a function $f$ with $g$ is the function $f * g$ defined by
$$
(f * g)(y) := \int_X f(z) \tau_y g(z) dvol(z)
$$
Fixing a basepoint $o \in X$, convolution with an $L^1$-function radial around $o$ gives rise to a bounded operator $T : L^p(X) \to L^p(X)$
for all $p \in [1, +\infty]$ satisfying the following properties (see section 2.4):

\medskip

\noindent (1) $T$ preserves the subspace of $L^p$-functions radial around $o$.

\medskip

\noindent (2) $T \tau_x \phi = \tau_x T \phi$ for all $\phi \in C^{\infty}_c(X)$ radial around $o$ and for all $x \in X$.

\medskip

\noindent (3) $T \phi * \psi = \phi * T \psi$ for all $\phi, \psi \in C^{\infty}_c(X)$ radial around $o$.

\medskip

This motivates the following definition: an $L^p$-multiplier is a bounded operator $T : L^p(X) \to L^p(X)$ which satisfies
properties (1)-(3) above. Examples of $L^p$-multipliers are given by convolution with radial $L^1$-functions, or more
generally convolution with radial complex measures of finite total variation (see section 2.4).

\medskip

The terminology "multiplier" is motivated by the following: for $p > 2$, if $T : L^p(X) \to L^p(X)$ is an $L^p$-multiplier,
then there exists a holomorphic function $m_T$ on the strip $S_p$, called the {\it symbol} of $T$, such that for all
$C^{\infty}_c$-functions $\phi$ radial around $o$, the spherical Fourier transform of $T \phi$ is given by
$$
\widehat{T\phi}^o(\lambda) = m_T(\lambda) \hat{\phi}^o(\lambda) \ , \lambda \in S_p
$$

Moreover if $T$ is not a scalar multiple of the identity, then the function $m_T$ is a nonconstant
holomorphic function. We can now state our main theorem:

\medskip

\begin{theorem} \label{mainthm} Let $X$ be a complete, simply connected, harmonic manifold with sectional curvature $K$
satisfying $K \leq -1$. Let $2 < p < \infty$ and let $T : L^p(X) \to L^p(X)$ be an $L^p$-multiplier
with symbol $m_T$ such that $T$ is not a scalar multiple of the identity.
Then for all $\lambda \in S_p$ such that $m_T(\lambda) \neq 0$, for any $\nu \in \C$ such that $|\nu| = |m_T(\lambda)|$
the dynamics of the operator $\frac{1}{\nu}T$ on $L^p(X)$ is topologically mixing with periodic points dense, in particular
the dynamics is chaotic in the sense of Devaney.
\end{theorem}

\medskip

A particular case of multipliers is given by the heat semigroup $e^{t\Delta}$ on $X$. For a simply connected harmonic
manifold, the heat kernel $H_t(x, y)$ is radial, i.e. there exists an $L^1$ function $h_t$ radial around $o$
such that $H_t(x, y) = (\tau_x h_t)(y)$ (see \cite{szabo}). The action of $e^{t\Delta}$ is thus given by
convolution with the radial $L^1$ function $h_t$, so $e^{t\Delta}$ is an $L^p$-multiplier for all $p \in [1, +\infty]$.
We determine the symbol of $e^{t\Delta}$ and then apply the previous theorem to obtain the following corollary:

\medskip

\begin{cor} \label{heat} Let $X$ be a complete, simply connected, harmonic manifold with sectional curvature $K$
satisfying $K \leq -1$, and let $2 < p < \infty, 1 < q < 2$ be such that $1/p + 1/q = 1$. There exists a constant $c_p = \frac{4\rho^2}{pq}$
such that the action of the shifted heat semigroup $(e^{ct}e^{t\Delta})_{t > 0}$ on
$L^p(X)$ is chaotic in the sense of Devaney for all $c \in \C$ with $\Re c > c_p$. In fact for any $t_0 > 0$, the operator
$e^{ct_0} e^{t_0 \Delta}$ on $L^p(X)$ is chaotic for all $c \in \C$ with $\Re c > c_p$.
\end{cor}

\medskip

In section 2 we recall some basic facts about eigenfunctions of the Laplacian, the Fourier transform, and convolution on
harmonic manifolds, show that convolution with a radial measure of finite variation is an
$L^p$-multiplier, and prove existence of the symbol of a multiplier. In section 3 we prove the main theorem.
We also prove the corollary by determining the symbol of the multiplier $e^{t\Delta}$.

\medskip

\section{Preliminaries}

\medskip

In this section we briefly recall the facts about the Fourier transform on harmonic manifolds which we will require.
For details the reader is referred to \cite{biswas6}.
Throughout, $X$ will denote a complete, simply connected harmonic $n$-manifold with sectional curvatures $K$
satisfying $K \leq -1$. We fix a basepoint $o \in X$.

\medskip

\subsection{CAT(-1) spaces and Busemann functions}

\medskip

In this case, $X$ is a CAT(-1) space, and we can define a {\it boundary at infinity} $\partial X$ of the space $X$,
defined as the set of equivalence classes of geodesic rays $\gamma : [0, \infty)$ in $X$, where two rays are equivalent if they
stay at bounded distance from each other. There is a natural topology on $\overline{X} := X \cup \partial X$ called the
{\it cone topology} for which $\overline{X}$ becomes a compactification of $X$ (for details on CAT($\kappa$) spaces we
refer to \cite{bridsonhaefliger}).

\medskip

Given a point $x \in X$, let $\lambda_x$ be normalized Lebesgue measure on the unit tangent sphere $T^1_x X$,
i.e. the unique probability measure on $T^1_x X$ invariant under the orthogonal group of the tangent space
$T_x X$. For $v \in T^1_x X$, let $\gamma_v : [0, \infty) \to X$ be the unique geodesic ray with initial velocity $v$.
Then we have a homeomorphism $p_x : T^1_x X \to \partial X, v \mapsto \gamma_v(\infty)$. The
visibility measure on $\partial X$ (with respect to the basepoint $x$)
is defined to be the push-forward $(p_x)_* \lambda_x$ of $\lambda_x$ under the map $p_x$; for
notational convenience, we will however denote the visibility measure on $\partial X$ by the same symbol $\lambda_x$.

\medskip

The {\it Busemann cocycle} $B : X \times X \times \partial X$ is defined by
$$
B(x, y, \xi) := \lim_{z \to \xi} (d(x, z) - d(y, z))
$$

Given a point $x \in X$ and a boundary point $\xi \in \partial X$, the {\it Busemann function at $\xi$ based at $x$} is defined by
$$
B_{\xi, x}(y) := B(y, x, \xi)
$$

The Busemann functions $B_{\xi, x}$ are $C^2$ convex functions, and their level sets are called {\it horospheres} based at $\xi$.

\medskip

\subsection{Radial and horospherical eigenfunctions of the Laplacian}

\medskip

Let $\Delta$ denote the Laplace-Beltrami operator of $X$, or Laplacian. As $X$ is harmonic, $X$ is also
{\it asymptotically harmonic}, i.e. all horospheres have constant mean curvature, so there is a constant
$h$ such that $\Delta B_{\xi, x} \equiv h$ for all $\xi \in \partial X, x \in X$. Since $X$ is negatively curved,
in fact $h > 0$. We let
$$
\rho := \frac{1}{2} h
$$

\medskip

A function $f$ on $X$ is called radial around a point $x \in X$ if $f$ is constant
on geodesic spheres centered at $x$. For any $x \in X$ and $\lambda \in \C$, there is a unique eigenfunction $\phi_{\lambda, x}$ of
$\Delta$ for the eigenvalue $-(\lambda^2 + \rho^2)$ which is radial around $x$ and satisfies $\phi_{\lambda, x}(x) = 1$. Moreover
for any fixed $y \in Y$, $\lambda \mapsto \phi_{\lambda, x}(y)$ is an entire function of $\lambda$.
The functions $\phi_{\lambda, x}$ are real-valued for $\lambda \in \R \cup i\R$, and bounded by $1$ for $|\Im \lambda| \leq \rho$.
Given $p > 2$, for all $\lambda$ in the strip $S_p := \{ |\Im \lambda| < (1 - 2/p)\rho \}$, the function $\phi_{\lambda, x}$ is
in $L^p(X)$.

\medskip

For any $x \in X, \xi \in \partial X$ and $\lambda \in \C$, the function $e^{(i\lambda - \rho)B_{\xi,x}}$ is an
eigenfunction of $\Delta$ for the eigenvalue $-(\lambda^2 + \rho^2)$. Note that this eigenfunction is constant
on horospheres based at $\xi$.

\medskip

\subsection{The spherical and Helgason Fourier transforms}

\medskip

Let $f \in L^1(X)$. Given a point $x \in X$, the {\it spherical Fourier transform of $f$ based at $x$} is the function $\hat{f}^x$ on $\R$
defined by pairing $f$ with the radial eigenfunctions $\phi_{\lambda,x}$:
$$
\hat{f}^x(\lambda) := \int_X f(y) \phi_{\lambda, x}(y) dvol(y) \ , \ \lambda \in \R
$$

\medskip

There exists a function $c$ on $\C - \{0\}$ satisfying, for some constants $C, K > 0$, the estimates
\begin{align*}
\frac{1}{C} |\lambda| & \leq |c(\lambda)|^{-1} \leq C |\lambda| , \quad \quad  0 < |\lambda|  \leq K \\
\frac{1}{C} |\lambda|^{(n-1)/2} & \leq |c(\lambda)|^{-1} \leq C |\lambda|^{(n-1)/2} , \quad |\lambda| \geq K \\
\end{align*}
such that the following inversion formula for the spherical Fourier transform from \cite{biswas6} holds:

\medskip

\begin{theorem} \label{inversion} Let $f \in C^{\infty}_c(X)$ be radial around $x$. Then
$$
f(y) = \int_{0}^{\infty} \hat{f}^x(\lambda) \phi_{\lambda, x}(y) |c(\lambda)|^{-2} d\lambda
$$
for all $y \in X$.
\end{theorem}

\medskip

Given $1 \leq q < 2$, if $p > 2$ is the conjugate exponent such that $1/p + 1/q = 1$, then using the fact that the functions
$\phi_{\lambda,o}$ are in $L^p(X)$ for $\lambda$ in the strip $S_p$, we have the following proposition from \cite{biswas6}:

\medskip

\begin{prop} \label{lpfourier} Let $1 \leq q < 2$ and $p > 2$ be such that $1/p + 1/q = 1$. Then for any $x \in X$ and $f \in L^q(X)$,
the spherical Fourier transform of $f$ based at $x$ is well-defined and extends to a holomorphic function on the strip $S_p$.
\end{prop}

\medskip

Let $f \in C^{\infty}_c(X)$. Given $x \in X$, the {\it Helgason Fourier transform of $f$ based at $x$} is the function
$\tilde{f}^x : \C \times \partial X \to \C$ defined by
$$
\tilde{f}^x(\lambda, \xi) := \int_X f(y) e^{(-i\lambda - \rho)B_{\xi,x}(y)} dvol(y) \ , \lambda \in \C, \xi \in \partial X
$$

\medskip

We have the following relation between the Helgason Fourier transforms
based at two different basepoints $o, x \in X$:

\medskip

\begin{equation} \label{fourierbasept}
\tilde{f}^x(\lambda, \xi) = e^{(i\lambda + \rho)B_{\xi, o}(x)} \tilde{f}^o(\lambda, \xi)
\end{equation}

\medskip

If $f$ is radial around the point $x$ then the Helgason Fourier transform redues to the spherical Fourier transform,
$$
\tilde{f}^x(\lambda, \xi) = \hat{f}^x(\lambda) \ , \lambda \in \C, \xi \in \partial X
$$

From \cite{biswas6} we have the following inversion formula for the Helgason Fourier transform:

\medskip

\begin{theorem} \label{inversion1} Let $x \in X$ and let $f \in C^{\infty}_c(X)$. Then
$$
f(y) = \int_{0}^{\infty} \int_{\partial X} \tilde{f}^x(\lambda, \xi) e^{(i\lambda - \rho)B_{\xi,x}(y)} d\lambda_x(\xi) |c(\lambda)|^{-2} d\lambda
$$
for all $y \in X$.
\end{theorem}

\medskip

\subsection{Convolution operators and $L^p$ multipliers}

\medskip

For a point $x \in X$, let $d_x$ denote the distance function from the point $x$, defined by $d_x(y) := d(x, y), y \in X$.

\medskip

Given a function $f$ on $X$ radial around a point $x$, let $u$ be a function on $[0, \infty)$ such that $f = u \circ d_x$.
Given a point $y$ in $X$, the $y$-translate of $f$ is the function $\tau_y f$ radial around
$y$ defined by $\tau_y f := u \circ d_y$. It follows from the fact that $X$ is harmonic that 
$||\tau_y f||_p = ||f||_p$ for all $p \in [1,+\infty]$. Moreover if
$f$ is also in $L^1$, then the spherical Fourier transforms satisfy
$$
\widehat{\tau_y f}^y(\lambda) = \hat{f}^x(\lambda)
$$
We note also from \cite{biswas6} that there is an even $C^{\infty}$ function on $\R$ which we denote by $\phi_{\lambda}$ such that
$\phi_{\lambda, x} = \phi_{\lambda} \circ d_x$. Thus the $x$-translate of the eigenfunction $\phi_{\lambda, o}$ radial
around $o$ is the eigenfunction $\phi_{\lambda, x}$ radial around $x$, $\tau_x \phi_{\lambda, o} = \phi_{\lambda, x}$.

\medskip

For simplicity, in the sequel, unless otherwise mentioned, by "radial function" we will mean a function which is radial around the
basepoint $o$. Likewise, by "spherical Fourier transform" we will mean the spherical Fourier transform based at $o$, unless
otherwise mentioned.

\medskip

Given $f,g \in L^1(X)$ with $g$ radial, the convolution of $f$ with $g$ is the function $f * g$ on $X$ defined by
$$
(f * g)(x) = \int_X f(y) \tau_x g(y) dvol(y)
$$
The integral above converges for a.e. $x$, and satisfies
$$
||f * g||_1 \leq ||f||_1 ||g||_1
$$
We note that if $f \in L^{\infty}(X)$ and $g \in L^1(X)$ with $g$ radial, then the integral defining $(f * g)(x)$ converges
for all $x$ and satisfies
$$
||f * g||_{\infty} \leq ||f||_{\infty} ||g||_1
$$
It follows by interpolation that for any $p \in [1, +\infty]$, convolution with a radial $L^1$ function $g$ defines a bounded
linear operator on $L^p(X)$ satisfying
$$
||f * g||_p \leq ||f||_p ||g||_1
$$
for all $f \in L^p(X)$.

\medskip

A standard argument using the above inequality and density of $C^{\infty}_c(X)$ in $L^p(X)$ gives that if $\{\phi_n\}$ is an
approximate identity, i.e. $\phi_n \geq 0, \int_X \phi_n dvol = 1$ and $\int_{B(o,r)} \phi_n dvol \to 1$ for any $r > 0$, then
for any $f \in L^p(X)$,
$$
||f * \phi_n - f||_p \to 0
$$
as $n \to \infty$.

\medskip

In \cite{biswas6} it is shown that for $\phi, \psi \in C^{\infty}_c(X)$ with $\psi$ radial, the Helgason Fourier transform
of the convolution $\phi * \psi$ satisfies
$$
\widetilde{\phi * \psi}^o(\lambda, \xi) = \tilde{\phi}^o(\lambda, \xi) \hat{\psi}^o(\lambda) \ , \lambda \in \C, \xi \in \partial X
$$

In particular, if both $\phi, \psi$ are radial, then
$$
\widehat{\phi * \psi}^o(\lambda) = \hat{\phi}^o(\lambda) \hat{\psi}^o(\lambda)
$$

We also have from \cite{biswas6} that the radial $L^1$ functions form a commutative Banach algebra under convolution. It follows,
using density of radial $C^{\infty}_c$-functions in radial $L^p$ functions,
that for a radial $L^1$ function $g$ the convolution operator $T_g :f \mapsto f * g$ on $L^p(X)$ preserves the subspace of
radial $L^p$ functions and satisfies, for all radial $C^{\infty}_c$-functions $\phi, \psi$,
$$
T_g \phi * \psi = \phi * T_g \psi
$$

In fact for any $x \in X$ the convolution operator $T_g$ preserves the subspace of $L^p$ functions radial around $x$. This is
a consequence of the following lemma:

\medskip

\begin{lemma} \label{transconv} Let $\phi, \psi$ be radial $C^{\infty}_c$-functions. Then for any $x \in X$,
$$
\tau_x \phi * \psi = \tau_x (\phi * \psi)
$$
\end{lemma}

\medskip

\noindent{\bf Proof:} We compute Helgason Fourier transforms:
\begin{align*}
\widetilde{\tau_x \phi * \psi}^o(\lambda, \xi) & = \widetilde{\tau_x \phi}^o(\lambda, \xi) \hat{\psi}^o(\lambda) \\
                                               & =  e^{-(i\lambda + \rho)B_{\xi, o}(x)} \widetilde{\tau_x \phi}^x(\lambda, \xi) \hat{\psi}^o(\lambda) \\
                                               & = e^{-(i\lambda + \rho)B_{\xi, o}(x)} \widehat{\tau_x \phi}^x(\lambda) \hat{\psi}^o(\lambda) \\
                                               & = e^{-(i\lambda + \rho)B_{\xi, o}(x)} \hat{\phi}^o(\lambda) \hat{\psi}^o(\lambda) \\
                                               & = e^{-(i\lambda + \rho)B_{\xi, o}(x)} \widehat{\phi * \psi}^o(\lambda) \\
                                               & = e^{-(i\lambda + \rho)B_{\xi, o}(x)} \widehat{\tau_x(\phi * \psi)}^x(\lambda) \\
                                               & = e^{-(i\lambda + \rho)B_{\xi, o}(x)} \widetilde{\tau_x(\phi * \psi)}^x(\lambda, \xi) \\
                                               & = \widetilde{\tau_x(\phi * \psi)}^o(\lambda, \xi) \\
\end{align*}

It follows from the Fourier inversion formula (Theorem \ref{inversion1}) that $\tau_x \phi * \psi = \tau_x(\phi * \psi)$. $\diamond$

\medskip

Now given $g$ a radial $L^1$ function and $\phi \in C^{\infty}_c(X)$, let $\{ \psi_n \}$ be a sequence of radial $C^{\infty}_c$-functions
converging to $g$ in $L^1$. Given $x \in X$, since $\phi$ and $\tau_x \phi$ are in $L^{\infty}$, it follows that
$\phi * \psi_n$ and $\tau_x \phi * \psi_n$ converge pointwise to $\phi * g$ and $\tau_x \phi * g$ respectively, so $\tau_x(\phi * \psi_n)$
converges pointwise to $\tau_x(\phi * g)$. Applying the previous Lemma, we obtain $\tau_x \phi * g = \tau_x(\phi * g)$. Thus
the convolution operator $T_g$ satisfies
$$
T_g \tau_x \phi = \tau_x T_g \phi
$$
for all radial $C^{\infty}_c$ functions $\phi$ and all $x \in X$.

\medskip

This leads us to the following definition:

\medskip

\begin{definition} ({\bf $L^p$-multipliers}) For $p \in [1, +\infty]$, an $L^p$-multiplier is a bounded
operator $T : L^p(X) \to L^p(X)$ such that:

\medskip

\noindent (1) $T$ preserves the subspace of radial $L^p$ functions.

\medskip

\noindent (2) For all radial $C^{\infty}_c$-functions $\phi, \psi$ we have
$$
T\phi * \psi = \phi * T\psi
$$

\medskip

\noindent (3) For all radial $C^{\infty}_c$-functions $\phi$ and all $x \in X$ we have
$$
T \tau_x \phi = \tau_x T \phi
$$
\end{definition}

\medskip

Thus convolution operators given by radial $L^1$ functions are $L^p$ multipliers for all $p \in [1, +\infty]$.
For more general examples of $L^p$-multipliers we
can consider convolution with radial complex measures $\mu$ of finite total variation, which is defined as follows:

\medskip

We say that a complex measure
$\mu$ on $X$ is radial around $o$ if there exists a complex measure $\tilde{\mu}$ on $[0, \infty)$ such that for any continuous
bounded function $f$ on $X$ we have
$$
\int_X f(x) d\mu(x) = \int_{0}^{\infty} \left( \int_{S(o,r)} f(y) d\lambda_{o,r}(y) \right) d\tilde{\mu}(r)
$$
where $S(o,r)$ denotes the geodesic sphere of radius $r$ around $o$ and $\lambda_{o,r}$ denotes the volume measure on $S(o,r)$ induced
from the metric on $X$. For $x \in X$, the $x$-translate of such a measure $\mu$ is the measure $\tau_x \mu$ radial around $x$ defined by
requiring that
$$
\int_X f(y) d\tau_x\mu(y) = \int_{0}^{\infty} \left( \int_{S(x,r)} f(y) d\lambda_{x,r}(y) \right) d\tilde{\mu}(r)
$$
for all continuous bounded functions $f$ on $X$ (where $S(x,r)$ is the geodesic sphere of radius $r$ around $x$ and $\lambda_{x,r}$ is the
volume measure on $S(x, r)$).

\medskip

For an $L^1$ function $f$ on $X$ and a radial complex measure $\mu$ on $X$ of finite variation, the convolution $f * \mu$ is the function
on $X$ defined by
$$
(f * \mu)(x) := \int_X f(y) d\tau_x \mu(y)
$$
We note that any $L^1$ function $g$ which is radial around $o$ gives a complex measure $\mu = g dvol$ which is radial around $o$ and satisfies
$||\mu|| = ||g||_1$ (where $||\mu||$ is the total variation norm of $\mu$), and $f * \mu = f * g$, so convolution with finite variation radial
measures generalizes convolution with $L^1$ radial functions.

\medskip

Given a finite variation radial measure $\mu$, we can approximate $\mu$ in the weak-* topology by measures $g_n dvol$ where $g_n$'s are
radial $L^1$ functions such that $||g_n||_1 \to ||\mu||$, then for any $f \in C^{\infty}_c(X)$ we have $f * g_n \to f * \mu$ pointwise,
and an application of Fatou's Lemma then leads to the inequality
$$
||f * \mu||_1 \leq ||f||_1 ||\mu||
$$
valid for all $f \in C^{\infty}_c(X)$ and all finite variation radial measures $\mu$. The inequality then continues to hold for all
$f \in L^1(X)$ by density of $C^{\infty}_c(X)$ in $L^1(X)$.

\medskip

Moreover for $f \in L^{\infty}(X)$ and $\mu$ a finite variation radial measure, it is straightforward to see that the integral
defining $f * \mu$ exists for all $x$ and satisfies
$$
||f * \mu||_{\infty} \leq ||f||_{\infty} ||\mu||
$$
Thus by interpolation for any $p \in [1, +\infty]$, convolution with a finite variation radial measure $\mu$ defines a bounded
operator on $L^p(X)$ satisfying
$$
||f * \mu||_p \leq ||f||_p ||\mu||
$$
for all $f \in L^p(X)$.

\medskip

\begin{prop} \label{lpmultiplier} Let $\mu$ be a radial complex measure of finite total variation. Then for any $p \in [1, +\infty]$,
the operator $T_{\mu} : f \mapsto f * \mu$ is an $L^p$ multiplier.
\end{prop}

\medskip

\noindent{\bf Proof:} Fix $p \in [1, \infty]$. Let $\{g_n\}$ be a sequence of radial $L^1$
functions such that $g_n dvol \to \mu$ in the weak-* topology and
such that $||g_n||_1 \to ||\mu||$. Then for any radial $C^{\infty}_c$-function $\phi$, the functions $\phi * g_n$ are
radial and converge to $\phi * \mu$ pointwise, so $\phi * \mu$ is radial. It follows that $T_{\mu}$ preserves
the subspace of radial $L^p$ functions.

\medskip

Let $\phi, \psi$ be radial $C^{\infty}_c$-functions. Then
$$
||\phi * g_n||_{\infty} \leq ||\phi||_{\infty} ||g_n||_1 \leq C ||\phi||_{\infty}
$$
for some constant $C > 0$, so for any $x \in X$ the functions $\phi * g_n$ are uniformly bounded on the support
of $\tau_x \psi$, and converge to $\phi * \mu$ pointwise, so it follows from dominated convergence that
$(\phi * g_n) * \psi(x) \to (\phi * \mu) * \psi(x)$ for all $x \in X$. A similar argument gives that
$\phi * (\psi * g_n)(x) \to \phi * (\psi * \mu)(x)$ for all $x \in X$. Since $(\phi * g_n) * \psi = \phi * (\psi * g_n)$ for all $n$,
it follows that $(\phi * \mu) * \psi = \phi * (\psi * \mu)$.

\medskip

Let $\phi$ be a radial $C^{\infty}_c$-function and let $x \in X$. Then $\phi * g_n$ and $\tau_x \phi * g_n$ converge to
$\phi * \mu$ and $\tau_x \phi * \mu$ respectively, so $\tau_x(\phi * g_n)$ converges pointwise to $\tau_x(\phi * \mu)$.
Since $\tau_x \phi * g_n = \tau_x(\phi * g_n)$ for all $n$, it follows that $\tau_x \phi * \mu = \tau_x(\phi * \mu)$.
$\diamond$

\medskip

Let $1 \leq q < 2$ and $p > 2$ such that $1/p + 1/q = 1$. Let $f$ be a radial $L^q$ function, then the spherical Fourier
transform $\hat{f}$ is holomorphic in the strip $S_p$, and it turns out that for any radial $C^{\infty}_c$-function
$\psi$, we have
$$
\widehat{f * \psi}(\lambda) = \hat{f}(\lambda) \hat{\psi}(\lambda) \ , \lambda \in S_p
$$
This can be seen as follows: let $\{ \phi_n \}$ be a sequence of radial $C^{\infty}_c$-functions converging to $f$
in $L^q(X)$, then since $\phi_{\lambda, o} \in L^p(X)$ for $\lambda \in S_p$, it follows from Holder's inequality
that $\hat{\phi_n}(\lambda) \to \hat{f}(\lambda)$ for $\lambda \in S_p$. Moreover, since $\psi \in L^1(X)$,
$\phi_n * \psi$ converges to $f * \psi$ in $L^q(X)$, so as before $\widehat{\phi_n * \psi}(\lambda) \to \widehat{f * \psi}(\lambda)$
for $\lambda \in S_p$. The desired equality follows by passing to the limit in the equality
$\widehat{\phi_n * \psi}(\lambda) = \hat{\phi_n}(\lambda) \hat{\psi}(\lambda)$.

\medskip

Other examples of $L^p$-multipliers can be obtained by using the {\it Kunze-Stein phenomenon} proved in \cite{biswas6}.
This asserts that if $1 \leq q < 2$, then there is a constant $C_q > 0$ such that for all $C^{\infty}_c$-functions
$f, g$ with $g$ radial, we have
$$
||f*g||_2 \leq C_q ||f||_2 ||g||_q.
$$
Combining this with the trivial estimate
$$
||f*g||_{\infty} \leq ||f||_{\infty} ||g||_1 ,
$$
it follows from interpolation that for any $p > 2$, if $1 \leq r < 2$ is such that $1/r < 1 + 1/p$, then there is a constant $C_p > 0$
such that
$$
||f*g||_p \leq C_p ||f||_p ||g||_r .
$$
The above inequality then implies that convolution with a radial $L^r$-function $g$ defines an $L^p$-multiplier $T_g : L^p(X) \to L^p(X)$.

\medskip

The following proposition justifies the use of the term "multiplier":

\medskip

\begin{prop} Let $1 \leq q < 2$ and $p > 2$ be such that $1/p + 1/q = 1$. Let $T : L^p(X) \to L^p(X)$ be an $L^p$ multiplier.
Then there exists a holomorphic function $m_T$ on the strip $S_p$ such that,
for any radial $C^{\infty}_c$-function $\phi$, we have $T\phi \in L^q(X)$, and
$$
\widehat{T \phi}(\lambda) = m_T(\lambda) \hat{\phi}(\lambda) \ , \lambda \in S_p
$$
\end{prop}

\medskip

\noindent{\bf Proof:} We first show that given a radial $C^{\infty}_c$ function $\phi$, $T \phi \in L^q(X)$.
For any radial $C^{\infty}_c$-function $\psi$, we have
\begin{align*}
\left|\int_X T \phi(x) \psi(x) dvol(x) \right| & = |T\phi * \psi(o)| \\
                                  & = |\phi * T\psi(o)| \\
                                  & = \left|\int_X \phi(x) T\psi(x) dvol(x)\right| \\
                                  & \leq ||\phi||_q ||T\psi||_p \\
                                  & \leq (||T|| ||\phi||_q) ||\psi||_p \\
\end{align*}
Since $T\phi$ is radial and the above inequality holds for all radial $C^{\infty}_c$-functions $\psi$, it follows that
$||T\phi||_q \leq ||T|| ||\phi||_q < +\infty$.

\medskip

Thus for any radial $C^{\infty}_c$-function $\phi$ which is not identically zero, $\widehat{T\phi}$ is a holomorphic
function in the strip $S_p$, and we can define a meromorphic function $m_{\phi}$ on $S_p$ by
$$
m_{\phi} := \frac{\widehat{T\phi}}{\hat{\phi}}
$$
If $\psi$ is another radial $C^{\infty}_c$-function which is not identically zero, then the equality
$T\phi * \psi = \phi * T\psi$ implies $\widehat{T\phi} \hat{\psi} = \hat{\phi} \widehat{T\psi}$ on $S_p$ and hence
$m_{\phi} = m_{\psi}$. Thus the meromorphic function $m_{\phi}$ is independent of the choice of $\phi$, and we may denote it
by $m_T$.

\medskip

It suffices to show that $m_T$ is in fact holomorphic in $S_p$. For this it is enough to show that given any $\lambda_0 \in S_p$,
there is a radial $C^{\infty}_c$-function $\phi$ such that $\hat{\phi}(\lambda_0) \neq 0$, since then $m_T = \widehat{T\phi}/\hat{\phi}$
will be holomorphic near $\lambda_0$. If $\hat{\phi}(\lambda_0) = 0$ for all radial $C^{\infty}_c$-functions $\phi$, then
$$
\int_X \phi(x) \phi_{\lambda_0, o}(x) dvol(x) = 0
$$
for all such $\phi$, and since $\phi_{\lambda_0, o}$ is radial this implies that $\phi_{\lambda_0, o} \equiv 0$, a contradiction.
Thus $m_T$ is holomorphic in $S_p$ and by definition satisfies $\widehat{T \phi} = m_T \hat{\phi}$ for all radial $C^{\infty}_c$-functions
$\phi$. $\diamond$

\medskip

\noindent{\bf Remark.} If for $1 \leq q < 2$ we have an $L^q$-multiplier $T$, then by definition $T \phi \in L^q$ for
$\phi$ a radial $C^{\infty}_c$-function, and then the proof of the above proposition applies to show that for
any $L^q$-multiplier $T$ there is a function $m_T$ holomorphic in the strip $S_p$ such that $\widehat{T\phi}(\lambda) = m_T(\lambda) \hat{\phi}(\lambda)$
for $\lambda \in S_p$ and $\phi$ a radial $C^{\infty}_c$-function. Thus the conclusion of the proposition holds
in fact for all $L^p$-multipliers with $p \neq 2$.

\medskip

We will call the holomorphic function $m_T$ given by the above proposition the {\it symbol} of the $L^p$-multiplier $T$.
Note that if $T$ is given by convolution with a radial $L^1$-function $g$, then the symbol $m_T$ equals the spherical
Fourier transform $\hat{g}^o$ of $g$, since $\widehat{\phi * g}^o = \hat{\phi}^o \hat{g}^o$ for all radial $C^{\infty}_c$-functions
$\phi$.

\medskip

\begin{prop} \label{eigenfunctions} Let $1 \leq q < 2$ and $p > 2$ be such that $1/p + 1/q = 1$. Let $T : L^p(X) \to L^p(X)$ be an $L^p$-multiplier. Then for all $\lambda \in S_p$ and $x \in X$, we have
$$
T \phi_{\lambda, x} = m_T(\lambda) \phi_{\lambda, x}
$$
\end{prop}

\medskip

\noindent{\bf Proof:} Let $\lambda \in S_p$ and let $\{ \phi_n \}$ be a sequence of radial $C^{\infty}_c$-functions converging to
$\phi_{\lambda, o}$ in $L^p(X)$. Then $T \phi_n$ converges to $T \phi_{\lambda, o}$ in $L^p(X)$. For any radial $C^{\infty}_c$-function
$\psi$, since $\psi \in L^q(X)$ it follows from Holder's inequality that
$$
\int_X T\phi_n(x) \psi(x) dvol(x) \to \int_X T\phi_{\lambda, o}(x) \psi(x) dvol(x)
$$
as $n \to \infty$. On the other hand, again using Holder's inequality and the fact that $\phi_n$ converges to $\phi_{\lambda, o}$
in $L^p(X)$, we have
\begin{align*}
\int_X T\phi_n (x) \psi(x) dvol(x) & = T\phi_n * \psi(o) \\
                                   & = \phi_n * T\psi(o) \\
                                   & = \int_X \phi_n(x) T\psi(x) dvol(x) \\
                                   & \to \int_X \phi_{\lambda, o} T\psi(x) dvol(x) \\
                                   & = \widehat{T\psi}(\lambda) \\
                                   & = m_T(\lambda) \hat{\psi}(\lambda) \\
                                   & = m_T(\lambda) \int_X \phi_{\lambda, o}(x) \psi(x) dvol(x) \\
\end{align*}
Thus
$$
\int_X T\phi_{\lambda, o}(x) \psi(x) dvol(x) = m_T(\lambda) \int_X \phi_{\lambda, o}(x) \psi(x) dvol(x)
$$
for all radial $C^{\infty}_c$-functions $\psi$, so it follows that $T\phi_{\lambda, o} = m_T(\lambda) \phi_{\lambda, o}$.

\medskip

Now given $x \in X$ and $\lambda \in S_p$, the functions $\tau_x \phi_n$ converge to $\phi_{\lambda, x}$ in $L^p(X)$,
and so
\begin{align*}
T \phi_{\lambda, x} & = \lim_{n \to \infty} T \tau_x \phi_n \\
                    & = \lim_{n \to \infty} \tau_x T \phi_n \\
                    & = \tau_x T \phi_{\lambda, o} \\
                    & = m_T(\lambda) \tau_x \phi_{\lambda, o} \\
                    & = m_T(\lambda) \phi_{\lambda, x} \\
\end{align*}
$\diamond$

\medskip

\section{Dynamics of $L^p$ multipliers}

\medskip

\subsection{General multipliers}

\medskip

We show in this section that the dynamics of appropriately scaled $L^p$-multipliers is chaotic in the sense
of Devaney if $2 < p < \infty$. The following lemma is the key to the results which follow:

\medskip

\begin{lemma} \label{lpdense} Let $1 < q < 2$ and $2 < p < \infty$ be such that $1/p + 1/q = 1$.
Let $K \subset S_p$ be a subset such that $K$ has a limit point in $S_p$. Then the subspace
$$
V_K := Span \{ \tau_x \phi_{\lambda, o} | x \in X, \lambda \in K \}
$$
is dense in $L^p(X)$.
\end{lemma}

\medskip

\noindent{\bf Proof:} It suffices to show that if $f \in L^q(X)$ is such that $\int_X f(y) \tau_x \phi_{\lambda, o}(y) dvol(y) = 0$
for all $x \in X, \lambda \in K$, then $f = 0$. Given such an $f \in L^q(X)$, the hypothesis on $f$ means that for any $x \in X$,
the spherical Fourier transform of $f$ based at $x$ vanishes on the set $K$. By Proposition \ref{lpfourier}, $\hat{f}^x$ is holomorphic in $S_p$ and
$K$ has a limit point in $S_p$, thus $\hat{f}^x$ vanishes identically in $S_p$, in particular on $\R$. Thus for all $x \in X$
and $\lambda \in \R$,
we have
$$
(f * \phi_{\lambda, o})(x) = \int_X f(y) \phi_{\lambda, x}(y) dvol(y) = \hat{f}^x(\lambda) = 0
$$
Let $\phi$ be a radial $C^{\infty}_c$-function, then by the Fourier inversion formula (Theorem \ref{inversion}) we have
$$
\phi(y) = \int_{0}^{\infty} \hat{\phi}(\lambda) \phi_{\lambda, o}(y) |c(\lambda)|^{-2} d\lambda
$$
for all $y \in X$, so it follows from Fubini's theorem that
$$
(f * \phi)(x) = \int_{0}^{\infty} (f * \phi_{\lambda, o})(x) \hat{\phi}(\lambda) |c(\lambda)|^{-2} d\lambda = 0
$$
for all $x \in X$. Thus $f * \phi = 0$ for all radial $C^{\infty}_c$-functions $\phi$. Now letting $\{ \phi_n \}$ be a
sequence of radial $C^{\infty}_c$-functions which forms an approximate identity, we have $f * \phi_n = 0$ for all $n$,
and $f * \phi_n$ converges to $f$ in $L^q(X)$, thus $f = 0$. $\diamond$

\medskip

We will also need the following lemma:

\medskip

\begin{lemma} \label{notscalar} Let $2 < p < \infty$ and let $T : L^p(X) \to L^p(X)$ be an $L^p$-multiplier. Suppose
$T$ is not a scalar multiple of the identity. Then the symbol $m_T$ is a nonconstant holomorphic function in the strip
$S_p$.
\end{lemma}

\medskip

\noindent{\bf Proof:} Suppose to the contrary that $m_T \equiv C$ for some constant $C \in \C$. By Proposition \ref{eigenfunctions} we
then have $T \phi_{\lambda, x} = C \phi_{\lambda, x}$ for all $\lambda \in S_p$ and $x \in X$. Thus $T = C Id$ on the subspace
$V = Span\{ \phi_{\lambda, x} | \lambda \in S_p, x \in X \}$, which is dense by the previous Lemma, hence $T = C Id$ on $L^p(X)$,
a contradiction. $\diamond$

\medskip

The main tool to prove that the dynamics of $L^p$ multipliers is chaotic is the following criterion of Godfrey-Shapiro
(\cite{erdmannmanguillot}, Theorem 3.1):

\medskip

\begin{theorem} \label{godshap} ({\bf Godfrey-Shapiro criterion}) Let $X$ be a separable Banach space and let $T : X \to X$
be a bounded operator. Suppose the subspaces $X^+, X^-$ defined by
\begin{align*}
X^+ & = Span\{ v \in X | Tv = \lambda v \ \hbox{for some } \lambda \in \C \ \hbox{such that } |\lambda| < 1 \} \\
X^- & = Span\{ v \in X | Tv = \lambda v \ \hbox{for some } \lambda \in \C \ \hbox{such that } |\lambda| > 1 \} \\
\end{align*}
are dense in $X$.
Then the dynamics of $T$ on $X$ is topologically mixing, i.e. for any two nonempty open sets $U, V \subset X$, there exists
$N \geq 1$ such that $T^n U \cap V \neq \emptyset$ for all $n \geq N$.
\end{theorem}

\medskip

We can now prove Theorem \ref{mainthm}:

\medskip

\noindent{\bf Proof of Theorem \ref{mainthm}:} Let $\lambda_0 \in S_p$ be such that $m_T(\lambda_0) \neq 0$, let
$\nu \in \C$ be such that $|\nu| = |m_T(\lambda)|$ and set $\alpha = m_T(\lambda)/\nu \in S^1$. Let $\mathbb{D}_0 = \{ z \in \C | |z| < 1 \}$
and $\mathbb{D}_{\infty} = \{ z \in \C |z| > 1 \}$. Let $U \subset S_p$ be an open neighbourhood of $\lambda_0$, then since
$\alpha \in S^1$ and by Lemma \ref{notscalar} $m_T$ is a nonconstant holomorphic function, there are nonempty open subsets
$U^+, U^- \subset U$ such that $\{ m_T(\lambda)/\nu | \lambda \in U^+ \} \subset \mathbb{D}_0$ and
$\{ m_T(\lambda)/\nu | \lambda \in U^- \} \subset \mathbb{D}_{\infty}$. By Proposition \ref{eigenfunctions}, for all
$\lambda \in U$ and $x \in X$, the function $\phi_{\lambda, x} \in L^p(X)$ is an eigenfunction of the operator $\frac{1}{\nu}T$
with eigenvalue $m_T(\lambda)/\nu$. By Lemma \ref{lpdense}, the subspaces $V^+ = \{ \phi_{\lambda, x} | \lambda \in U^+, x \in X \}$
and $V^- = \{ \phi_{\lambda, x} | \lambda \in U^-, x \in X \}$ are dense in $L^p(X)$. It follows from the Godfrey-Shapiro
criterion that the dynamics of $\frac{1}{\nu}T$ is topologically mixing.

\medskip

It remains to show that the periodic points of $\frac{1}{\nu}T$ are dense in $L^p(X)$. Since $m_T$ is a nonconstant holomorphic function
and $m_T(\lambda_0)/\nu \in S^1$, we can choose sequences $\{\lambda_n\} \subset U$ and $\{ p_n/q_n \} \subset \mathbb{Q}$
such that $m_T(\lambda_n)/\nu = e^{2\pi i p_n/q_n}$ and $\lambda_n \to \lambda_0$ as $n \to \infty$. Then by Lemma \ref{lpdense},
the subspace $V = Span\{ \phi_{\lambda_n, x} | x \in X, n \geq 1 \}$ is dense in $L^p(X)$. It thus suffices to show that each
element of $V$ is a periodic point of $\frac{1}{\nu}T$. Any element $\phi \in V$ can be written as
$\phi = \sum_{j = 1}^N a_j \phi_{\lambda_j, x_j}$ for some $N \geq 1, a_1, \dots, a_N \in \C$ and $x_1, \dots, x_N \in X$.
Since $\phi_{\lambda_j, x_j}$ is an eigenvector of $\frac{1}{\nu}T$ with eigenvalue $e^{2\pi i p_j/q_j}$, letting $q = \prod_{j = 1}^N q_j$
it follows that $(\frac{1}{\nu}T)^q \phi_{\lambda_j, x_j} = \phi_{\lambda_j, x_j}$ for all $j$,
thus $(\frac{1}{\nu}T)^q \phi = \phi$ and
$\phi$ is a periodic point of $\frac{1}{\nu}T$. $\diamond$

\medskip

\subsection{The heat semigroup}

\medskip

We recall some basic facts about the heat semigroup and heat kernel on a complete Riemannian manifold $X$.
Denote by $\Delta_X = div \ grad$ the Laplacian acting on $C^{\infty}_c(X) \subset L^2(X)$, then this is an
essentially self-adjoint operator, and so its closure $\Delta_{X, 2}$ is a self-adjoint operator on
$L^2(X)$. Since $\Delta_{X, 2}$ is negative, it generates a semigroup $e^{t\Delta_{X, 2}}$ on
$L^2(X)$ by the spectral theorem for unbounded self-adjoint operators. The operators $e^{t\Delta_{X, 2}}$
are positive, leave $L^1(X) \cap L^{\infty}(X) \subset L^2(X)$ invariant, and may be extended
to a positive contraction semigroup $e^{t\Delta_{X, p}}$ on $L^p(X)$ for any $p \in [1, +\infty]$, which
is strongly continuous for $p \in [1, +\infty)$ (\cite{davies}). In the sequel we will write
simply $e^{t\Delta}$ for the semigroup $e^{t\Delta_{X, p}}$ on $L^p(X)$.
From \cite{strichartz} we have the
following:

\medskip

There exists a $C^{\infty}$ function $H_t(x, y)$ on $\R^+ \times X \times X$, the {\it heat kernel}, such that
for all $t > 0$ and $x \in X$ the function $H_t(x, .)$ is positive and in $L^p$ for all $p \in [1, +\infty]$, and
for all $f \in L^p(X)$,
$$
e^{t\Delta} f(x) = \int_X f(y) H_t(x,y) dvol(y)
$$
and
$$
\frac{\partial}{\partial t} e^{t\Delta} f(x) = \Delta e^{t\Delta} f(x)
$$
Moreover, it is shown in \cite{szabo} that for a $X$ a simply connected harmonic manifold, the heat kernel is radial,
i.e. there exists a function $h_t$ radial around the basepoint $o$ such that $H_t(x, y) = (\tau_x h_t)(y)$. Thus the action
of the heat semigroup on $L^p(X)$ is given in our case by convolution with the radial $L^1$ function $h_t$,
$$
e^{t\Delta} f = f * h_t
$$
for all $f \in L^p(X)$, so $e^{t\Delta}$ is an $L^p$-multiplier for all $p \in [1, +\infty]$.
The symbol of the multiplier $e^{t\Delta}$ is given by the following proposition:

\medskip

\begin{prop} \label{heatfourier} For any $t > 0$, the spherical Fourier transform of the heat kernel is given by
$$
\hat{h_t}^o(\lambda) = e^{-t(\lambda^2 + \rho^2)} \ , \lambda \in S_{\infty}
$$
\end{prop}

\medskip

\noindent{\bf Proof:} Let $p \in (2, \infty)$ and let $\lambda \in S_p$. Then $\phi_{\lambda, o} \in L^p(X)$, and
using the fact that the operators $\Delta, e^{t\Delta}$ on $L^p(X)$ commute and $\Delta \phi_{\lambda, o} = -(\lambda^2 + \rho^2) \phi_{\lambda, o}$,
we have
\begin{align*}
\frac{\partial}{\partial t} e^{t\Delta} \phi_{\lambda, o} & = \Delta e^{t\Delta} \phi_{\lambda, o} \\
                                                          & = e^{t\Delta} \Delta \phi_{\lambda, o} \\
                                                          & = -(\lambda^2 + \rho^2) e^{t\Delta} \phi_{\lambda, o} \\
\end{align*}

%
%
%
%
Thus $t \mapsto e^{t\Delta} \phi_{\lambda, o} \in L^p(X)$ satisfies the first order linear ODE
$$
\frac{\partial}{\partial t} e^{t\Delta} \phi_{\lambda, o} = -(\lambda^2 + \rho^2) e^{t\Delta} \phi_{\lambda, o}
$$
and $e^{t\Delta} \phi_{\lambda, o} \to \phi_{\lambda, o}$ in $L^p(X)$ as $t \to 0$, hence
$$
e^{t\Delta} \phi_{\lambda, o} = e^{-t(\lambda^2 + \rho^2)} \phi_{\lambda, o}
$$
for all $t > 0$. Evaluating both sides above at the point $o$ gives
\begin{align*}
\hat{h_t}^o(\lambda) & = \int_X \phi_{\lambda, o}(x) h_t(x) dvol(x) \\
                     & = e^{t\Delta} \phi_{\lambda, o}(o) \\
                     & = e^{-t(\lambda^2 + \rho^2)} \phi_{\lambda, o}(o) \\
                     & = e^{-t(\lambda^2 + \rho^2)} \\
\end{align*}
$\diamond$

\medskip

We can now prove the result on the chaotic dynamics of shifted heat semigroups:

\medskip

\noindent{\bf Proof of Corollary \ref{heat}:} Given $2 < p < \infty$ and $1 < q < 2$ such that $1/p + 1/q = 1$, let
$c_p = 4\rho^2/(pq)$. Let $c \in \C$ be such that $\Re c > c_p$, and let $t_0 > 0$. Let $T = e^{t_0 \Delta}$ and $\nu = e^{-ct_0}$.
By Proposition \ref{heatfourier} above, the symbol of $T$ is given by $m_T(\lambda) = e^{-t_0 (\lambda^2 + \rho^2)}$.
In order to show that the operator $e^{ct_0} e^{t_0 \Delta} = \frac{1}{\nu} T$ is chaotic, it suffices by Theorem \ref{mainthm}
to show that there exists $\lambda \in S_p$ such that $|\nu| = |m_T(\lambda)|$.

\medskip

Letting $\lambda = s + it \in S_p$, the equality $|\nu| = |m_T(\lambda)|$ is equivalent to
$$
s^2 - t^2 + \rho^2 = \Re c
$$
Let $t$ be such that $t = (1 - 2/p)\rho - \epsilon$ where $\epsilon > 0$ is small, then we have
\begin{align*}
\Re c + t^2 - \rho^2 & = (\Re c - c_p) + c_p + ((1 - 2/p)^2 - 1)\rho^2 + O(\epsilon) \\
                     & = (\Re c - c_p) + (4(1/p)(1 - 1/p) - 4/p + 4/p^2)\rho^2 + O(\epsilon) \\
                     & = (\Re c - c_p) + O(\epsilon) \\
                     & > 0 \\
\end{align*}
for $\epsilon$ small enough since $\Re c - c_p > 0$. Thus we can choose $t$ with $0 < t < (1 - 2/p)\rho$ such that $\Re c + t^2 - \rho^2 > 0$,
so we can then choose $s \in \R$ such that $s^2 = \Re c + t^2 - \rho^2$, or $s^2 - t^2 + \rho^2 = \Re c$, as required. $\diamond$

\bibliography{moeb}
\bibliographystyle{alpha}

\end{document}